\documentclass{amsart}
\usepackage{amssymb,amsfonts,amsmath}
\usepackage[all]{xy}
\usepackage[shortlabels]{enumitem}
\usepackage{hyperref, color}		
\newtheorem{theorem}{Theorem}[section]

\newtheorem{corollary}{Corollary}[section]
\theoremstyle{definition}

\theoremstyle{remark}
\newtheorem{remark}{Remark}
\theoremstyle{claim}

\theoremstyle{theoremA}

\theoremstyle{theoremB}
  
\newcommand{\real}{\mathbb{R}}

\newcommand{\la}{\lambda}
\newcommand{\La}{\Lambda}

\newcommand{\om}{\omega}

\newcommand{\na}{\nabla}

\newcommand{\ga}{\gamma}
\newcommand{\Ga}{\Gamma}

\newcommand{\ka}{\kappa}
\newcommand{\De}{\Delta}

\newcommand{\varep}{\varepsilon}

\newcommand{\dd}{\,\textrm{d}M}
\newcommand{\ti}{\tilde}
\newcommand{\vol}{\text{\rm vol}}

\newcommand{\lan}{\left\langle}
\newcommand{\ran}{\right\rangle}

\newcommand{\m}{\mathcal}
\newcommand{\tr}{\mathrm{tr}\,}
\newcommand{\ric}{\mathrm{Ric}\,}
\newcommand{\spp}{\mathrm{supp}}
\newcommand{\II}{I\!I}
\begin{document}


\title[$L^2$-harmonic 1-forms on submanifolds]
{$L^2$-harmonic 1-forms on submanifolds \\ with finite total curvature}
\author{M. P. Cavalcante}
\thanks{The  authors were partially supported by CNPq-Brazil}

\author{H. Mirandola}

\author{F. Vit\'orio}

\date{}

\maketitle
\thispagestyle{empty}
\begin{abstract} Let $x:M^m\to \bar M$, $m\geq 3$, be an isometric immersion of a 
complete noncompact manifold $M$  in a complete simply-connected 
manifold $\bar M$ with sectional curvature satisfying $-k^2\leq K_{\bar M}\leq 0$, 
for some constant $k$. 
Assume that the immersion has finite total curvature in the sense that the
traceless second fundamental form has finite $L^m$-norm. If $K_{\bar M}\not \equiv 0$, assume further that the 
first eigenvalue of the Laplacian of $M$ is bounded from below by a suitable constant.  
We prove that the space of the $L^2$ harmonic $1$-forms on $M$ has finite dimension.  
Moreover there exists a constant $\La>0$, explicitly computed, such that if the total curvature 
is bounded from above by $\La$ then there is no nontrivial $L^2$-harmonic $1$-forms on $M$.
\end{abstract}
\section{Introduction}

Let $x:M^m\to \bar M$ be an isometric immersion of an $m$-dimensional   
manifold $M$  in a Riemannian manifold $\bar M$. Let us denote  $\II$ 
the second  fundamental form and $H=\frac{1}{m}\tr(\II)$ the mean curvature 
vector field of the immersion $x$. The traceless second fundamental form $\Phi$ is defined by
\begin{equation*}
\Phi(X,Y)=\II(X,Y) - \lan X,Y\ran H,
\end{equation*}
for all vector fields $X,Y$ on $M$, where  $\lan\,,\ran$ is the metric of $M$.  A simple computation shows that
\begin{equation*}
|\Phi |^2= |\II |^2-m|H|^2.
\end{equation*}
In particular, $|\Phi|\equiv 0$ if and only  if the immersion $x$ is totally umbilical. 
We say that the immersion $x$ has {\it finite total curvature} if the $L^m$-norm of the 
traceless second fundamental form is finite (see \cite{BdCS} and  references therein), that is,
\begin{equation*}
\|\Phi\|_{L^m(M)}=\left(\int_M |\Phi|^m  \dd\right)^{\frac{1}{m}} < +\infty,
\end{equation*}
where $\dd$ stands for the volume element of $M$.

Topological and metric properties of complete submanifolds with finite total curvature 
have been a field of active 
research, since the work of Gauss. For instance, let $M^2$ be a complete surface isometrically immersed in an  
Euclidean space $\real^n$ with finite total curvature. By celebrated results of Huber \cite{h},  Osserman \cite{o} 
and Chern-Osserman \cite{co}, 
it is known that if the immersion is minimal then it is proper, $M$ is homeomorphic to a compact surface 
$\m M$ punctured at a finitely many points and the Gauss map  extends continuously to all points of $\m M$. 
See also the papers of White \cite{wh} and M\"uller-\u Sver\'ak \cite{ms} for the non-minimal case. 
In higher dimension, the description of the topology is more involved and 
there exist a lot of interesting papers related with this subject (see, for instance, \cite{a}, \cite{dccs}, \cite{do}, \cite{bc}, \cite{bmj}).

In the present paper, we are interested in the study of cohomological aspects of 
noncompact submanifolds with finite total curvature. More specifically, 
assume that  $M$ is a complete noncompact manifold and 
consider the space of the  $L^2$-harmonic $1$-forms on $M$
\begin{equation*}
\m H^1(M) = \Big\{\om \ \bigm|  \int_M \om \wedge*\om = \int_M |\om|^2 \dd
 <\infty \mbox{ and } \,\textrm{d}\om=\textrm{d}\! *\!\om = 0 \Big\}.
\end{equation*}
It is well known that the space $\m H^1(M)$ is isomorphic to the first-reduced $L^2$-cohomology 
group of $M$ (see section 3.b of \cite{c}). Moreover, the dimension of $\m H^1(M)$ gives an 
upper bound to the number of non-parabolic ends of $M$.  In fact, if $\m H^0_D(M)$ 
denotes the space of the harmonic functions on $M$ with finite Dirichlet integral then 
$u\in \m H^0_D(M)$ if and only if its differential exterior $du\in \m H^1(M)$. Thus
\begin{equation*}\label{nonpar-ends}
\dim \m H_D^0(M)\leq \dim \m H^1(M)+1.
\end{equation*}
On the other hand, an important result of Li-Tam (Theorem 2.1 of \cite{lt}) 
states that the number $\mathbf{e}(M)$ of 
nonparabolic ends of a complete manifold $M$ satisfies $\mathbf{e}(M)\leq \dim \m H^0_D(M)$. 

In \cite{c}, Theorem 3.5, Carron  proved that if $M^m$, $m\geq3$, is a complete noncompact 
submanifold of $\real^n$  with finite total curvature and finite total mean curvature 
(i.e.,  the $L^m$-norm of the mean  curvature vector is finite) then 
each  space of reduced $L^2$-cohomology on $M$ has finite dimension.
Under same conditions, and using techniques of  harmonic function as in \cite{lt},  
Fu and Xu  also proved that $\m H^1(M)$ is finite dimensional (see Theorem 1.1 of \cite{fx}).
Since the ends of complete noncompact submanifolds in $\real^n$ with 
finite total mean curvature are  nonparabolic they actually conclude that $M$ 
must have finitely many ends.  

Our first result is an improvement and a generalization of Carron and Fu-Xu theorems.
We recall that a Riemannian manifold is called a Hadamard manifold if it is complete, 
simply-connected and has nonpositive sectional curvature. 

\begin{theorem}\label{t1} 
Let $x:M^m\to \bar M$, $m\geq 3$, be an isometric immersion of a complete noncompact 
manifold $M$ in  a Hadamard manifold $\bar M$ with sectional 
curvature satisfying $-k^2\leq K_{\bar M}\leq 0$, for some constant $k$.  
In the case $K_{\bar M}\not \equiv  0$, assume further that the first eigenvalue of the 
Laplace-Beltrami  operator of $M$ satisfies
\begin{equation*}
\la_1(M)>\frac{(m-1)^2}{m}(k^2- \liminf_{r(p)\to \infty}|H(p)|^2),
\end{equation*}
where $r$ stands for the distance in $M$ from a fixed point. If $x$ has finite total curvature then 
the space  $\m H^1(M)$ has finite dimension.
\end{theorem}
\medskip
An interesting result of Anderson \cite{a2} shows that, for all $m\geq 3$ and 
$\ka>m-2$, there exists a complete simply-connected manifold $M_\ka^m$ 
with sectional curvature satisfying $-\ka^2\leq K\leq -1$ and such that 
$\dim \m H^1(M^m_\ka) = \infty$.
Since $\bar M=M^m_\ka\times \real$ is a Hadamard manifold satisfying 
$-\ka^2\leq K_{\bar M}\leq 0$
and $M^m_\ka$ is a totally geodesic submanifold of $\bar M$ we conclude 
that the hypothesis on the first eigenvalue in Theorem \ref{t1} is necessary.
%

\medskip

In the same paper \cite{c} Carron proved a gap theorem on the dimension of $\m H^1(M)$.
More precisely, Carron proved that there exists a constant $\varep(m)$ 
such that if $\| \II\|\leq \varep(m)$, then all spaces of $L^2$-harmonic forms are trivial. 
In the next result we prove a gap theorem for immersion with small $\|\Phi\|_{L^m(M)}$.
This result is also a generalization of \cite{fl}.

\begin{theorem}\label{t2} 
Let $x:M^m\to \bar M$, $m\geq 3$, be an isometric immersion of a complete noncompact 
manifold $M$ in  a Hadamard manifold $\bar M$ with sectional 
curvature satisfying $-k^2\leq K_{\bar M}\leq 0$, for some constant $k$.  
In the case $K_{\bar M}\not \equiv  0$, assume further that the first eigenvalue of the 
Laplace-Beltrami  operator of $M$ satisfies
\begin{equation*}
\la_1(M)>\frac{(m-1)^2}{m}(k^2 - \inf |H|^2).
\end{equation*}
There exists a positive constant $\La$ 
such that if $\|\Phi\|_{L^m(M)}<\La$ then there is no nontrivial $L^2$-harmonic $1$-form on $M$.
Furthermore, if $k=0$ then $\La$ depends only on $m$; 
otherwise, $\La$ depends only on $m$, k,  $\la_1(M)$ and $\inf|H|$.
\end{theorem}
It is a natural question to ask about the best constant $\La$ in Theorem \ref{t2}. 
In the next result we present an explicit value for $\La$ depending on the case. 
\begin{theorem}\label{t3} 
Let $S=S(m,2)$ be the constant of Sobolev's inequality derived from \cite{hs}{\footnote{Let $x:M^m\to \bar M$ be an isometric immersion of a complete manifold $M$ in a Hadamard manifold $\bar M$. 
By a simple computation, it shows that for all $1\leq p<m$ the following inequality holds:
\begin{equation}\label{soblv}
\Big(\int_M h^{\frac{pm}{m-p}}dV\Big)^\frac{m-p}{m} \leq S(m,p) \int_M \Big(|\na h|^ p+ (h|H|)^p\Big)dV,
\end{equation}
for all nonnegative $C^ 1$-function $h:M^ m\to \real$ with compact support, where $S(m,p)^{\frac{1}{p}}=\frac{2p(m-1)}{m-p}c(m)$ and $c(m)$ is the positive constant, depending only on $m$, that appears in Theorem 2.2 of \cite{hs}. 
}}. Then
the constant $\La$ as in Theorem \ref{t2} can be given explicitly as follows:
\begin{enumerate}[(i)]
\item\label{i} 
If $k=0$ and $H=0$ then 
$$
\La=\frac{m}{(m-1)\sqrt{S}}.
$$ 
\item\label{ii} 
If $k=0$ and $H$ is arbitrary then 
$$
\La=\dfrac{\dfrac{m}{(m-1)\sqrt{S}}}{\sqrt{1+\dfrac{(m-1)(m-2)^2}{4(m^2-3m+1)}}}.
$$  
\item\label{iii} 
If $k\neq 0$ and $H=0$ then 
$$
\La=\frac{m}{(m-1)\sqrt{S}}\sqrt{1-\frac{(m-1)^2 k^2}{m \la_1(M)}}.
$$
\item\label{iv} 
If $k\neq 0$ and $\inf|H|>\Big(1-\dfrac{m}{(m-1)^2}\Big)^{-\frac{1}{2}}|k|$ then 
$$
\La=\dfrac{\frac{m}{(m-1)\sqrt{S}}}
{\sqrt{1+\dfrac{(m-2)^2(m-1)^2}{4m((m^2-3m+1)-\dfrac{(m-1)^2 k^2}{\inf|H|^2})}}}.
$$
\item\label{v} 
If $k\neq 0$, $H$ is arbitrary and $\la_1(M)>\dfrac{(m-1)^2k^2}{m}$ then 
$$
\La=\dfrac{\dfrac{m}{(m-1)\sqrt{S}}\sqrt{1-\dfrac{(m-1)^2 k^2}{m \la_1(M)}}}
{\sqrt{1+\dfrac{(m-1)(m-2)^2}{4(m^2-3m+1+\dfrac{(m-1)^2k^2}{\la_1(M)})}}}.
$$
\item\label{vi} 
If $k\neq 0$, $\inf|H|\leq \Big(1-\frac{m}{(m-1)^2}\Big)^{-\frac{1}{2}}|k|$ and 
$\la_1\leq \dfrac{(m-1)^2k^2}{m}$ then  
$$
\La=\frac{m}{(m-1)\sqrt{S}}
\sqrt{\dfrac{1-\m B\dfrac{(m-2)^2(m-1)}{4m}\inf|H|^2}{(1+\dfrac{1}{\m B})(\la_1(M)+\inf|H|^2)}},
$$ 
where 
$\m B=-1+\sqrt{1+\frac{4m\m A}{(m-1)(m-2)^2\inf|H|^2}}$ and  
$\m A=\la_1-\frac{(m-1)^2}{m}(k^2-\inf|H|^2)$.
\end{enumerate}
\end{theorem}
\begin{remark}
Item \ref{i} was previously obtained by Seo in \cite{seo} as an improvement of a 
a result due to Ni (\cite{ni}). Item \ref{iii} was also obtained by Seo \cite{seo2} 
in the case that $\bar M$ is the hyperbolic space $\mathbb{H}^n$.  Item \ref{iii} and Item \ref{i} are covered by Theorem 3.2 and Corollary 3.12, respectively, in H. P. Fu: \emph{ Studies on geometry and analysis of submanifolds}, Ph. D. dissertation, Zhejiang University, 2007.
Note that if $k\neq 0$ and $\frac{(m-1)^2}{m}(k^2-\inf \|H\|^2)< \la_1\leq \frac{(m-1)^2k^2}{m}$ 
then $\inf|H|>0$, 
and hence $\m B$ is well defined.
\end{remark}
As a direct consequence of Theorems \ref{t2} and \ref{t3} we obtain the following results: 
\begin{corollary} 
Let $x:M^m\to \real^n$, $m\geq 3$, be an isometric immersion of a complete noncompact 
manifold $M^m$  in the Euclidean space $\real^n$. 
If the total curvature of $x$ satisfies 
$$
\|\Phi\|_{L^m(M)}<\dfrac{\dfrac{m}{(m-1)\sqrt{S}}}{\sqrt{1+\dfrac{(m-1)(m-2)^2}{4(m^2-3m+1)}}}
$$
then there is no nontrivial $L^2$-harmonic $1$-forms on $M$.  
\end{corollary}

\begin{corollary} Let $x:M^m\to \bar M$, $m\ge 3$,  be an isometric immersion of a 
complete noncompact manifold $M$ in a Hadamard manifold $\bar M$ 
with sectional curvature satisfying $-k^2\leq K_{\bar M}\leq 0$, for some constant $k\neq 0$. 
Assume that $\la_1(M)>\frac{(m-1)^2}{m}k^2$. If the total curvature of $x$ satisfies 
$$
\|\Phi\|_{L^m(M)}<
\dfrac{\dfrac{m}{(m-1)\sqrt{S}}\sqrt{1-\dfrac{(m-1)^2 k^2}
{m \la_1(M)}}}{\sqrt{1+\dfrac{(m-1)(m-2)^2}{4(m^2-3m+1+\dfrac{(m-1)^2k^2}{\la_1(M)})}}}
$$ 
then there is no nontrivial $L^2$-harmonic $1$-forms on $M$.
\end{corollary}

This paper is organized as follows. In section \ref{desigualdade} we obtain a general
inequality for $1$-forms on $M$ involving the geometry of immersion. 
In section  \ref{s3} we use such inequality to prove Theorem \ref{t2}. 
In section \ref{s4} we prove Theorem \ref{t3}. 
Finally, in section \ref{s5} we use again the inequality of section \ref{desigualdade}
to  give a proof of Theorem \ref{t1}.


\section{An integral inequality for $L^2$-harmonic $1$-forms}\label{desigualdade}

 Let $x:M^m\to \bar M$, $m\ge 3$,  be an isometric immersion of a 
complete noncompact manifold $M$ in a Hadamard manifold $\bar M$ 
with sectional curvature satisfying $-k^2\leq K_{\bar M}\leq 0$, for some constant $k$. 
In this section, we obtain an integral inequality for $L^2$-harmonic $1$-forms on 
$M$ involving the geometry of $x$.

Given $\omega \in \m H^1(M)$, we  recall the refined Kato's inequality (see, for instance, 
Lemma 3.1 of \cite{w}):
\begin{equation*}
|\nabla |\om||^2\leq \frac{m-1}{m}|\nabla \om|^2.
\end{equation*}
A direct computation yields
\begin{equation*} \label{straight}
\Delta |\omega|^2\leq 2\bigl( |\om| \De |\om|+\frac{m-1}{m}|\nabla \om|^2 \bigr).
\end{equation*}
Using the Bochner's Formula \cite{b} (see also Lemma 3.2  of \cite{li}), we obtain
\begin{equation} \label{kato}
|\om|\De|\om| \geq \frac{1}{m-1}  |\nabla |\om||^2+\ric_M(\om^\#,\om^\#),
\end{equation}
where $\om^\#$ is the dual vector field of $\om$.
Under our hypothesis on the sectional curvature of $\bar M$ 
we can estimate the Ricci curvature of $M$ by using Proposition 2 of \cite{sx}:
\begin{eqnarray}\label{shio}
\ric_M(\om^\#,\om^\#) &\geq&(m-1)(|H|^2-k^2)|\om|^2 \\&&\,-\frac{(m-1)}{m}|\Phi|^2|\om|^2 -
\frac{(m-2)\sqrt{m(m-1)}}{m}|H| |\Phi||\om|^2 \nonumber.
\end{eqnarray}
Using (\ref{kato}) and (\ref{shio}), we get
\begin{eqnarray} \label{bks}
|\om| \De |\om| &\geq &
\frac{1}{m-1}|\nabla |\om||^2 -\frac{m-1}{m}|\Phi|^2|\omega|^2+
(m-1)\left(|H|^2-k^2\right)|\omega|^2\\
&&-\frac{(m-2)\sqrt{m(m-1)}}{m}|H| |\Phi||\omega|^2.\nonumber
\end{eqnarray}
Let $\eta\in C^\infty_0(M)$ be a smooth function on $M$ with compact support. 
We multiply both sides of (\ref{bks}) by $\eta^2$ and integrate by parts. 
For the sake of simplicity, henceforth we will omit the element volume 
in the integrals. So, we obtain
\begin{eqnarray}\label{fx1}
 0&\leq& -2\int_{M}\eta |\om| \lan \nabla \eta, \nabla |\om| \ran  - 
\frac{m}{m-1}\int_{M} \eta^2 |\nabla |\om||^2 \nonumber\\ 
&&+ \frac{(m-2)\sqrt{m(m-1)}}{m} \int_{M} \eta^2 |H||\Phi||\om|^2
+\frac{m-1}{m}\int_{M} \eta^2|\Phi|^2|\om|^2\\
&&+ (m-1)\int_{M} \eta^2\left(k^2 - |H|^2 \right)|\om|^2.\nonumber
\end{eqnarray}
For each $a>0$, we apply the Cauchy-Schwarz Inequality in (\ref{fx1}) to obtain
\begin{eqnarray}\label{fx2}
 0&\leq& -2\int_{M}\eta |\om| \lan \nabla \eta, \nabla |\om| \ran  - 
\frac{m}{m-1}\int_{M} \eta^2 |\nabla |\om||^2\nonumber\\ 
&&+ \int_{M} \left((m-1)k^2 +
\bigg(-(m-1) + \frac{a(m-2)\sqrt{m(m-1)}}{2m}\bigg)|H|^2\right)  
\eta^2 |\om|^2\\
&&+\left( \frac{(m-2)\sqrt{m(m-1)}}{2a m}+
\frac{m-1}{m} \right) \int_{M} \eta^2|\Phi|^2|\om|^2.\nonumber
\end{eqnarray}
On the other hand, since $m\geq 3$, we use Holder, Hoffman-Spruck 
\cite{hs} and Cauchy-Schwarz Inequalities to get
\begin{eqnarray}\label{H-S}
\int_{M}\eta^2|\Phi|^2 |\om|^2 
&\leq& 
\phi(\eta)\left(\int_{M} (\eta |\om|)^{\frac{2m}{m-2}}\right)^{\frac{m-2}{m}}\nonumber\\
&\leq& S\phi(\eta) \int_{M} \Big(|\nabla(\eta |\om|)|^2 + \eta^2 |\om|^2 |H|^2\Big)\\
&\leq& S\phi(\eta) \int_{M} \Big((1+\frac{1}{b})|\om|^2 |\nabla \eta|^2 + 
(1+b)\eta^2 |\nabla|\om||^2 \Big)\nonumber\\&&
+ S\phi(\eta) \int_{M} \eta^2 |\om|^2 |H|^2,\nonumber
\end{eqnarray} 
for all $b>0$, where $\phi(\eta)=\left(\int_{\spp(\eta)} |\Phi|^m \right)^{\frac{2}{m}}$ and 
$S=S(m)>0$ is the constant in the Hoffman-Spruck Inequality. 
Thus, using (\ref{fx2}) and (\ref{H-S}), we have
\begin{eqnarray}\label{FX}
 0&\leq& 
 -2\int_{M}\eta |\om| \lan \nabla \eta, \nabla |\om| \ran  - 
 \frac{m}{m-1}\int_{M} \eta^2 |\nabla |\om||^2
+ A(m,a) \int_{M} |H|^2 \eta^2 |\om|^2\nonumber\\
&&+SB(m,a)\phi(\eta) \int_{M} \Big((1+\frac{1}{b})|\om|^2 |\nabla \eta|^2 + 
(1+b)\eta^2 |\nabla|\om||^2 \Big)\\
&&+SB(m,a) \phi(\eta) \int_{M} \eta^2 |\om|^2 |H|^2 + (m-1)k^2\int_M \eta^2|\om|^2,\nonumber
\end{eqnarray}
where $A(m,a)$ and $B(m,a)$ are given by
\begin{equation}\label{A-B}
\begin{array}{rcl}
A(m,a)&=& -(m-1)+\dfrac{a(m-2)\sqrt{m(m-1)}}{2m}\\ \\
B(m,a)&=& \dfrac{(m-2)\sqrt{m(m-1)}}{2a m}+\dfrac{m-1}{m}.
\end{array}
\end{equation}
We use the Cauchy-Schwarz Inequality again to get
\begin{equation}\label{cauchy_be}
2\left|\int_{M}\eta|\om|\lan\nabla\eta,\nabla |\om|\ran\right| 
\leq c \int_{M} \eta^2 |\nabla |w||^2 +\frac{1}{c}\int_{M}|\om|^2|\nabla \eta|^2,
\end{equation}
for all $c>0$. Using (\ref{FX}) and (\ref{cauchy_be}), we obtain the following integral inequality:
\begin{equation}\label{Equation}
C \int_M \eta^2 |\nabla |\om||^2 + D \int_M |H|^2\eta^2|\om|^2 
\leq E\int_M |\om|^2 |\nabla \eta|^2 + (m-1)k^2\int_{M} \eta^2 |\om|^2,
\end{equation}
where 
\begin{equation}\label{CDE}
\begin{array}{rcl}
-C&=&-C(m,a,b,c,\eta)=c+ (1+b)SB(m,a)\phi(\eta) - \frac{m}{m-1},\\ \\
-D&=&-D(m,a,\eta)=A(m,a) +SB(m,a)\phi(\eta),\\ \\
E&=&E(m,a,b,c,\eta)= \frac{1}{c}+(1+\frac{1}{b})SB(m,a)\phi(\eta).
\end{array}
\end{equation}

\section{Proof of Theorem \ref{t2}}\label{s3}
In this section, we will prove the existence of a positive constant $\La$ 
such that if $\|\Phi\|_{L^m(M)}<\La$ then $\m H^1(M)=\{0\}$. 
Choose $0<d<\frac{1}{2}$, $a=a(d)>0$ and $\La=\La(d)>0$ satisfying: 
\begin{equation}\label{estimate-ga}
\left\{
\begin{array}{l}
d+(m-1)(1+d)d<\frac{m}{m-1},\\ \\
\dfrac{a(m-2)\sqrt{m(m-1)}}{2m}<(m-1)d,\\ \\
SB(m,a)\La^2<(m-1)d.
\end{array}\right.
\end{equation}
Now we set
\begin{equation}\label{ti C}
\begin{array}{rcl}
-\bar C&=&
-\bar C(m,\La,a,b,c)=
c+ (1+b)SB(m,a)\La^2 - \frac{m}{m-1}, \mbox{ and }\\ \\ 
-\bar D&=&-\bar D(m,\La,a)=A(m,a)+SB(m,a)\La^2.
\end{array}
\end{equation}

Using (\ref{A-B}) and choosing  $0<c<d$ and $0<b<d$ we get:
\begin{equation*} 
\begin{array}{rcl}
\bar C &>& \frac{m}{m-1}-d-(m-1)(1+d)d>0,\\ \\
\bar D &>& (m-1)(1-2d)>0.
\end{array}
\end{equation*}

Assume that the total curvature of $x$ satisfies $\|\Phi\|_{L^m(M)}\leq \La$. 
Plugging the above choices in (\ref{Equation}) we obtain
\begin{equation}\label{ti_C}
\bar C \int_M \eta^2 |\nabla |\om||^2 + \bar D \int_M |H|^2\eta^2|\om|^2 
\leq \bar E\int_M |\om|^2 |\nabla \eta|^2 +(m-1)k^2\int_{M} \eta^2 |\om|^2,
\end{equation}
where $\bar E= \frac{1}{c}+(1+\frac{1}{b})SB(m,a)\La^2.$

In particular, if $k=0$, we obtain
\begin{equation}\label{kizero}
\bar C \int_M \eta^2 |\nabla |\om||^2 + \bar D \int_M |H|^2\eta^2|\om|^2 
\leq \bar E\int_M |\om|^2 |\nabla \eta|^2.
\end{equation}
We will see later that this inequality is sufficient to prove our result in the case $k=0$. 
In this case, we also note that $\La=\La(\ga)$ depends only on $m$. 
In order  to deal in the case $k\neq 0$ we need to introduce a new ingredient. 
We recall that the first eigenvalue $\la_1=\la_1(M)$ of the Laplacian of $M$ satisfies
\begin{equation}\label{lap}
\la_1 \int_M \varphi^2 \leq \int_M |\na \varphi|^2,
\end{equation}
for all $\varphi\in C^\infty_0(M)$. Applying  (\ref{lap}) with $\varphi=\eta|\om|$ and once more
using Cauchy-Schwarz inequality  we get for all $e>0$
\begin{equation*}
\la_1 \int_M \eta^2|\om|^2
\leq (1+e)\int_M \eta^2 |\na|\om||^2 + (1+\frac{1}{e})\int_M |\om|^2|\na \eta|^2,
\end{equation*}
which implies that 
\begin{equation}\label{xi}
\frac{\bar C\la_1}{1+e}\int_M \eta^2|\om|^2 
\leq \bar C \int_M \eta^2 |\na|\om||^2 + \frac{\bar C}{e}\int_M |\om|^2|\na \eta|^2.
\end{equation}
Thus, using (\ref{ti_C}) and (\ref{xi}), we obtain
\begin{eqnarray*}
\frac{\bar C\la_1}{1+e}\int_M \eta^2|\om|^2 
&\leq& \Big((m-1)k^2 - \bar D\inf|H|^2\Big)\int_{M} \eta^2 |\om|^2 \nonumber\\
&& + \Big(\bar E+\frac{\bar C}{e}\Big) \int_{M} |\om|^2 |\nabla \eta|^2.
\end{eqnarray*}
Note that 
\begin{eqnarray}\label{c_lambda1}
\frac{\bar C\la_1}{1+e}-(m-1)k^2 + \bar D\inf|H|^2
&\geq& \Big(\frac{m}{m-1}-d-(m-1)(1+d)d\Big)\frac{\la_1}{1+e}\nonumber\\
&& - (m-1)\Big(k^2 - (1-2d)\inf|H|^2\Big).
\end{eqnarray}
Thus, if $\la_1>\frac{(m-1)^2}{m}(k^2-\inf|H|^2)$ then we can choose $d$ and $e$, 
sufficiently small and depending on $m$, $k^2$, $\la_1$ and $\inf|H|^2$, so that 
$\frac{\bar C\la_1}{1+e}-(m-1)k^2 + \bar D\inf|H|^2>0$. Hence we get
\begin{equation}\label{ineq_c>0}
\int_M |\om|^2\eta^2 \leq \bar F \int_{M} |\om|^2 |\nabla \eta|^2,
\end{equation}
for some constant $\bar F>0$. In this case ($k\neq 0$), 
the constant $\La$ depends on $m$, $k^2$, $\la_1$ and $\inf|H|^2$.

\medskip

For each $r>0$, let $B_r$ denote the geodesic ball of radius $r$ on $M$ centered at 
some fixed point and let $\eta_r\in C^\infty_0(M)$ be a smooth function such that
\begin{equation*}\left\{
\begin{array}{lll}
0\leq \eta_r\leq 1& \mbox{ in }M,\\ 
\eta_r = 1& \mbox{ in }B_r,\\
|\nabla\eta_r|\leq 2r^{-1} & \mbox{ in }M,\\ 
\spp(\eta_r)=B_{2r}.
\end{array}\right.
\end{equation*}

If $k=0$ we use (\ref{kizero}) with $\eta_r$ to obtain
\begin{equation*}
\bar C\int_{B_r} |\nabla |\om||^2 + \bar D\int_{B_r} |H|^2 |\om|^2 
\leq \bar E\frac{4}{r^2} \int_{M} |\om|^2.
\end{equation*}
Taking $r\to \infty$ we get $|H||\om|=|\nabla|\om||=0$. 
Thus $|\om|$ is constant.  If $\om$ is not identically zero then $H=0$. 
In this case, since $\bar M$ is a Hadamard manifold it is well known that $M$ has 
infinite volume,  which is a contradiction, since $\int_M|\om|^2 <\infty$.
If $k\neq 0$ then using (\ref{ineq_c>0}) with $\eta_r$ we get
\begin{equation*}
\int_{B_r} |\om|^2 \leq \bar F \frac{4}{r^2} \int_{M} |\om|^2,
\end{equation*}
Taking $r\to \infty$ we have $\om=0$ and it finishes the proof.


\section{Proof of Theorem \ref{t3}}\label{s4} 
The key of the proof of Theorem \ref{t2} is the fact that there exists a constant 
$\La>0$ such that if $\|\Phi\|_{L^m(M)}<\La$ then one of the inequalities (\ref{kizero}) or
 (\ref{ineq_c>0}) holds. 
 We also proved that $\La$ depends only on $m$ in the case $k=0$ and 
 $\La$ depends on $m$, $k^2$, $\la_1$ and $\inf|H|$ in the case $k\neq 0$. 
 The goal of the present section is to give explicit estimates to $\La$.

\medskip

Using (\ref{A-B}) and (\ref{ti C}) we have
\begin{equation*}
\begin{array}{rcl}
-\bar C(m,\La,a,0,0)&=& S\La^2B(m,a) - \dfrac{m}{m-1}\\
&=&S\La^2 \Big(1+ \dfrac{(m-2)}{2a}\sqrt{\dfrac{m}{m-1}}\Big)\frac{m-1}{m} - \dfrac{m}{m-1} \,\textrm{ and } \\ \\
-\bar D(m,\La,a)&=& A(m,a)+S\La^2 B(m,a)\\
&=&\dfrac{m-1}{m}\Big(-m+\dfrac{a(m-2)}{2}\sqrt{\frac{m}{m-1}}+ 
S\La^2\Big(1+\dfrac{m-2}{2a}\sqrt{\frac{m}{m-1}}\Big)\Big).
\end{array}
\end{equation*}
Thus, using (\ref{ti C}) and the continuity of $\bar C$, it follows that there exist 
$\La>0$ and $a>0$ satisfying $\bar C=\bar C(m,\La,a,b,c)>0$, 
for some $b>0$ and $c>0$, sufficiently small, if and only if 
\begin{equation}\label{est-f1}
S\La^2 < f_1(a):=\frac{\frac{m^2}{(m-1)^2}}{1+\frac{m-2}{2a}\sqrt{\frac{m}{m-1}}}.
\end{equation}
Note that the function $f_1:(0,\infty)\to \real$ is increasing and $\sup f_1=\frac{m^2}{(m-1)^2}$.
Assume that \[\|\Phi\|_{L^m(M)}<\La_1:=\frac{m}{(m-1)\sqrt{S}}\] 
and take $\La>0$ so that $\|\Phi\|_{L^m(M)}<\La<\La_1$.
Since $S\La^2<\sup f_1$ and $\phi(\eta)<\La^2$, there exists $a_1=a_1(m,\La)>0$ 
such that, for any $a>a_1$, there exist $b>0$ and $c>0$, sufficiently small, 
satisfying $C>\bar C>0$ (see (\ref{CDE}) and (\ref{ti C})).  
Thus, if $k=0$ and $H=0$, we obtain from (\ref{Equation}) that inequality (\ref{kizero}) holds.
Thus, Item \ref{i} is proved.
\medskip

Similarly, there exist $\La>0$ and $a>0$ such that $\bar D=\bar D(m,\La,a)>0$ 
if and only if
\begin{equation}\label{est-f2}
S\La^2 < f_2(a):=\frac {m-\frac{a(m-2)}{2}\sqrt{\frac{m}{m-1}}}{1+\frac{m-2}{2a}\sqrt{\frac{m}{m-1}}}.
\end{equation}
The function $f_2:(0,\infty)\to \real$ is concave and 
$\max f_2 = f_2(\sqrt{\frac{m}{m-1}})=\frac{m}{m-1}> f_1(\sqrt{\frac{m}{m-1}})=\frac{2m}{(m-1)^2}$. 
Thus, the maximum value $\La_2>0$ that satisfies $\bar C>0$ and $\bar D>0$, 
for any $0<\La<\La_2$, and for some $a>0$, $b>0$ and $c>0$, is obtained in
the intersection point $a_{12}$ of the graphs of $f_1$ and $f_2$. Namely,
$$
a_{12}=\frac{2m(m^2-3m+1)}{(m-2)(m-1)^2\sqrt{\frac{m}{m-1}}}.
$$  
Hence we set
\begin{equation*}
\La_2:=\sqrt{\frac{f_2(a_{12})}{S}}=
\frac{m}{(m-1)\sqrt{S}}\Big(1+\frac{(m-2)^2(m-1)}{4(m^2-3m+1)}\Big)^{-\frac{1}{2}}.
\end{equation*}
If $\|\Phi\|_{L^m(M)}<\La<\La_2$, then, taking $a=a_{12}$, we have $C>\bar C>0$ and $D>\bar D>0$, 
for suitable constants $b>0$ and $d>0$. 
Therefore, for  $k=0$ and $H$ is arbitrary, inequality (\ref{kizero}) holds. Thus, Item \ref{ii} is proved.

\medskip
Now, we deal with the case $k\neq 0$ and $H=0$.  
Since we are assuming $\la_1>\frac{(m-1)^2}{m}(k^2-\inf|H|^2)$
we immediately obtain that $\la_1>0$ and $0<\frac{(m-1)^2}{m}\frac{k^2}{\la_1}<1$.
Using (\ref{ti C}) and the continuity of $\bar C$, we see that there exist 
$\La>0$ and $a>0$ satisfying $\bar C>0$ and 
$\frac{\bar C\la_1}{1+e}-(m-1)k^2>0$, 
for some $e>0$, $b>0$ and $c>0$, sufficiently small, if and only if $\La<\La_1$ and 
\begin{equation}\label{f3}
S\La^2<f_3(a):=
\frac{m^2}{(m-1)^2}\dfrac{1-\frac{(m-1)^2}{m}\dfrac{k^2}{\la_1}}
{1+\dfrac{(m-2)}{2a}\sqrt{\dfrac{m}{m-1}}}.
\end{equation}
The function $f_3:(0,\infty)\to \real$ is increasing and 
$\sup f_3=\frac{m^2}{(m-1)^2}\big(1-\frac{(m-1)^2}{m}\frac{k^2}{\la_1}\big)$. 

Set 
$\La_3:=\frac{m}{(m-1)\sqrt{S}}\Big(1-\frac{(m-1)^2}{m}\frac{k^2}
{\la_1}\Big)^{\frac{1}{2}}
$ 
and suppose that $\|\Phi\|_{L^m(M)}<\La<\La_3$.
Since $S\La^2<\sup f_3<S\La_3^2<S\La_1^2$ and $\phi(\eta)<\La^2$, 
there exists $a_3=a_3(m,k,\la_1,\La)>0$ such that, for any $a>a_3$, 
there exist $e>0$, $b>0$ and $c>0$, sufficiently small, satisfying $C>\bar C>0$ 
and $\frac{\bar C \la_1}{1+e}-(m-1)k^2>0$. 
Thus, in the case $k\neq 0$ and $H=0$ inequality (\ref{ineq_c>0}) holds and 
hence Item \ref{iii} follows. 

\medskip

Now, we assume that $H$ is arbitrary. Using (\ref{ti C}) 
we have that there exist $\La>0$ 
and $a>0$ such that $\bar C>0$, $\bar D>0$ and 
$\frac{\bar C\la_1}{1+e} - (m-1)k^2 + \bar D\inf|H|^2>0$, 
for some $e>0$, $b>0$ and $c>0$, sufficiently small, if and only if $\La<\La_2$ and 
\begin{equation*}\label{f4}
S\La^2<f_4(a):=
\frac{\frac{m^2}{(m-1)^2}
\big(\la_1-\frac{(m-1)^2}{m}\big(k^2-\inf|H|^2\big)\big)-
\frac{a(m-2)}{2}\inf|H|^2\sqrt{\frac{m}{m-1}}}
{\big(\la_1+\inf|H|^2\big)\big(1+\frac{(m-2)}{2a}\sqrt{\frac{m}{m-1}}\big)}.
\end{equation*}
By a simple computation, we can rewrite $f_4$ by the following ways:
\begin{eqnarray}
f_4&=&\frac{\inf|H|^2}{\la_1+\inf|H|^2}f_2 + 
\frac{1}{\la_1+\inf|H|^2}\Big(\la_1-\frac{(m-1)^2}{m}k^2\Big)f_1\label{f4 infH}\nonumber\\ 
&=& f_1 + \frac{\inf|H|^2}{\la_1+\inf|H|^2}(f_2-f_1) - 
\frac{\frac{(m-1)^2}{m}k^2}{\la_1+\inf|H|^2}f_1\label{f4-f1}\\ 
&=& f_2 + \frac{\la_1}{\la_1+\inf|H|^2}(f_1-f_2) - 
\frac{\frac{(m-1)^2}{m}k^2}{\la_1+\inf|H|^2}f_1.   \label{f4-f2} 
\end{eqnarray}
Using (\ref{f4-f1}) we obtain that $f_4$ intersects $f_1$ if and only if 
$\inf|H|^2>(1-\frac{m}{(m-1)^2})^{-1}k^2$. 
In this case, the intersection point of $f_1$ and $f_4$ is 
$$
a_{14}=\frac{2m\Big(m^2-3m+1 -
\dfrac{(m-1)^2k^2}{\inf|H|^2}\Big)}{(m-2)(m-1)^2\sqrt{\dfrac{m}{m-1}}}<a_{12}.
$$ 
Using that $f_1<f_2$ in $(0,a_{12})$, we obtain from (\ref{f4-f2}) that
 $f_4<f_2$ in $(0,a_{12})$.  
Thus, the maximum value $\La_4>0$ such that $\bar C>0$, $\bar D>0$ and 
 $\frac{\bar C\la_1}{1+e} - (m-1)k^2 + \bar D\inf|H|^2>0$, for all $0<\La<\La_4$ 
 and for some $a>0$, $b>0$, $c>0$ and $e>0$ is obtained considering 
\begin{equation}
S\La_4^2=f_4(a_{14})= 
\frac{m^2}{(m-1)^2}\Bigg(1+\frac{(m-2)^2(m-1)^2}{4m((m^2-3m+1)-
\dfrac{(m-1)^2k^2}{\inf|H|^2})}\Bigg)^{-1}.
\end{equation} 
Thus, assume that $\inf|H|^2>(1-\frac{m}{(m-1)^2})^{-1}k^2$ and $\|\Phi\|_{L^m(M)}<\La_4$.
Take $\La>0$ so that $\|\Phi\|_{L^m(M)}<\La<\La_4$ and let $a=a_{14}$.  
We obtain that $C>\bar C>0$, $D>\bar D>0$ and $\frac{\bar C\la_1}{1+e} - (m-1)k^2 + 
\bar D\inf|H|^2>0$, for some $e>0$, $b>0$ and $c>0$, sufficiently small. 
This implies that
the inequality (\ref{ineq_c>0}) holds. Item \ref{iv} is proved.

\medskip

Now, assume that $ \inf|H|^2\leq (1-\frac{m}{(m-1)^2})^{-1}k^2$. 
This implies that $f_4<f_1$, since $f_4$ does not intersect $f_1$ and 
$\lim_{a\to \infty}f_4(a)<\lim_{a\to \infty}f_1(a)$. 
Using (\ref{f4-f2}) we obtain that $f_4$ intersects $f_2$ if and only if $\la_1>0$. 
In this case, the intersection point of $f_4$ and $f_2$ is 
$$
a_{24}:=\frac{2m}{(m-2)\sqrt{\frac{m}{m-1}}}
\Big(1-\frac{m}{(m-1)^2}+\frac{k^2}{\la_1}\Big).
$$ 
Note also that $f_2(a)>0$ if and only if $0<a<\frac{2m}{(m-2)\sqrt{\frac{m}{m-1}}}$. 
This implies that $f_4(a_{24})>0$ if and only if $\la_1>\frac{(m-1)^2}{m}k^2$. 
Thus, assume further that $\la_1>\frac{(m-1)^2}{m}k^2$ and 
$\|\Phi\|_{L^m(M)}<\La_5:=\sqrt{\frac{f_4(a_{24})}{S}}$. 
Take $\La>0$ so that $\|\Phi\|_{L^m(M)}<\La<\La_4$ and let $a=a_{24}$. 
Using that $\phi(\eta)<\La^2$ and $S\La^2<f_2(a_{24})=f_4(a_{24})<f_1(a_{24})$  
we obtain that $C>\bar C>0$, $D>\bar D>0$ and 
$\frac{\bar C\la_1}{1+e} - (m-1)c^2 + \bar D\inf|H|^2>0$, 
for some $e>0$, $b>0$ and $c>0$, sufficiently small. 
This implies that inequality (\ref{ineq_c>0}) holds. It proves Item \ref{v}. 

To finish the proof of Theorem \ref{t3} we assume that 
$k\neq 0$, $\inf|H|^2\leq (1-\frac{m}{(m-1)^2})^{-1}k^2$ and 
$\la_1\leq \frac{(m-1)^2}{m}k^2$. 
This implies that $f_4<f_1$ and $f_4<f_2$. 
Furthermore, it holds that $\inf|H|>0$, since 
$\frac{(m-1)^2}{m}(k^2-\inf|H|^2)<\la_1\leq \frac{(m-1)^2}{m}k^2$. 
Thus, the maximum value $\La_6$ such that $\bar C>0$, $\bar D>0$ and 
$\frac{\bar C\la_1}{1+e} - (m-1)k^2 + \bar D\inf|H|^2>0$, for some suitable constants 
$a>0$, $b>0$, $c>0$ and $e>0$, is obtained considering $S\La_6^2=\sup f_4$. 
Note also that $\sup f_4=f_4(a_4)$, where  
$$
a_4=\frac{(m-2)}{2}\sqrt{\frac{m}{m-1}}\m B,
$$ 
$\m B=-1+\sqrt{1+\frac{4m\m A}{(m-1)(m-2)^2\inf|H|^2}}$ and 
$\m A=\la_1-\frac{(m-1)^2}{m}(k^2-\inf|H|^2)$. 
Thus, if 
$$
\|\Phi\|_{L^m}<\sqrt{\frac{f_4(a_4)}{S}}=
\frac{m}{(m-1)\sqrt{S}}\sqrt{\frac{1-\m B\frac{(m-2)^2(m-1)}{4m}
\inf|H|^2}{(1+\frac{1}{\m B})(\la_1+\inf|H|^2)}},
$$ 
then, using similar arguments to the previous case and taking $a=a_4$, 
we  obtain suitable positive constants $b$, $c$ and $e$ such that
 $C>\bar C>0$, $D>\bar D>0$ and $\frac{\bar C\la_1}{1+e} - (m-1)k^2 + \bar D\inf|H|^2>0$. 
 This implies that inequality (\ref{ineq_c>0}) holds. Item \ref{vi} is proved.


\section{Proof of Theorem \ref{t1}}\label{s5}
Consider $0<d<\frac{1}{2}$, $a=a(d)>0$ and $\La=\La(a)>0$ as given in 
(\ref{estimate-ga}). Assume that the immersion $x$ has finite total curvature.
Fix $r_0>0$ so that 
\begin{equation} \|\Phi\|_{L^m(M -B_{r_0})}<\sqrt{\La}.
\end{equation} 

Let $\eta=\eta_{r_0}\in C^\infty_0(M)$ be any smooth function with compact
support satisfying $\spp(\eta)\subset M-B_{r_0}$. 
Since $\phi(\eta)=\|\Phi\|^2_{L^m(\spp(\eta))}<\La$, we can proceed similarly as 
in section \ref{s3} (inequality (\ref{ti_C})) to obtain 
\begin{equation}\label{ineq_th-2}
C\int_M \eta^2 |\nabla |\om||^2 + \ti D\int_M |H|^2 \eta^2 |\om|^2 
\leq E \int_{M} |\om|^2 |\nabla \eta|^2 + (m-1)k^2\int_M \eta^2|\om|^2,
\end{equation}
for all $\om\in \m H^1(M)$, where  $\ti D=(m-1)(1-2d)$.

To deal in the case $k\neq 0$, we assume further that 
$$\la_1>\frac{(m-1)^2}{m}(k^2- \liminf_{r(p)\to \infty}|H(p)|^2),$$
where $r$ stands for the distance in $M$ from a fixed point.
It is easy to see that we can also consider $r_0>0$ sufficiently large satisfying
\begin{equation*}
\la_1>\frac{(m-1)^2}{m}(k^2-\inf_{M-B_{r_0}}|H|^2).
\end{equation*}
As in section \ref{s3} we obtain
\begin{equation}\label{ineq_c>0-2}
\int_M |\om|^2 \eta^2 \leq \ti F \int_{M} |\om|^2 |\nabla \eta|^2,
\end{equation}
for some constant $\ti F>0$.

It follows from Cauchy-Schwarz and Hoffman-Spruck \cite{hs} inequalities the following
\begin{eqnarray}\label{fx6}
S^{-1}\left(\int_{M} (\eta |\om|)^{\frac{2m}{m-2}}\right)^{\frac{m-2}{m}} 
&\leq&\int_{M} |\nabla(\eta |\om|)|^2 + \int_{M}\eta^2 |\om|^2|H|^2\nonumber\\
&\leq& (1+{s})\int_{M}\eta^2 |\nabla |\om||^2 + 
\big(1+\frac{1}{s}\big)\int_{M}|\om|^2 |\nabla\eta|^2\\
&&+\int_{M} \eta^2 |\om|^2|H|^2,\nonumber
\end{eqnarray}  
for any $s>0$. 
Using (\ref{ineq_th-2}) and (\ref{fx6}) we obtain
\begin{eqnarray}\label{sobol}
S^{-1}\left(\int_{M} (\eta |\om|)^{\frac{2m}{m-2}}\right)^{\frac{m-2}{m}}
&\leq& \big((1+s)C^{-1}\ti D+1\big)\int_M |H|^2 \eta^2|\om|^2 \nonumber\\
&& + \big(1+\frac{1}{s}\big)\int_M |\om|^2|\na \eta|^2 \\
&&+ (1+s)C^{-1}(m-1)k^2\int_M \eta^2|\om|^2,\nonumber
\end{eqnarray}
for all $s>0$. 
Since $m\geq 3$, we can choose $d$ and $s$ sufficiently small so that  
\begin{equation}\label{nu}
(1+s)C^{-1}\ti D = (1+s) \Big(\frac{m}{m-1}-d-(m-1)(1+d)d\Big)^{-1} (m-1)(2d-1) < -1.
\end{equation}
Thus, using (\ref{ineq_c>0-2}) (in the case $k\neq 0$), 
(\ref{sobol}) and (\ref{nu}) we obtain the following inequality
\begin{equation}\label{estimate}
S^{-1}\left(\int_{M} (\eta |\om|)^{\frac{2m}{m-2}}\right)^{\frac{m-2}{m}} 
\leq \ti A \int_M |\om|^2|\na \eta|^2,
\end{equation}
for some constant $\ti A=\ti A(m)>0$, for all $\om\in \m H^1(M)$.

>From now on, the proof follows standard technics (for instance, 
 \cite{lw} after equation 2.7, 
 or \cite{fx} after equation (18),  
 or \cite{cwx} after equation 3.14) 
 and uses a  {Moser iteration} argument and Lemma 11 of \cite{l2}. 
We include the proof here for the  sake of  completeness. 

Take $r>r_0+1$ and let $\eta$ be a smooth function satisfying the following conditions:

\begin{equation*}\left\{
\begin{array}{lll}
\eta=0& \mbox{ in }B_{r_0}\cup (M-B_{2r}),\\
\eta = 1& \mbox{ in }B_r-B_{r_0+1},\\
|\nabla\eta_r|\leq c_1 & \mbox{ in }B_{r_0+1}-B_{r_0},\\
|\nabla\eta_r|\leq c_1r^{-1} & \mbox{ in } B_{2r}-B_r,
\end{array}\right.
\end{equation*}
for some positive constant $c_1$.
Since $\spp(\eta)\subset M-B_{r_0}$, it follows from  (\ref{estimate}) the following
\begin{equation*}
\left(\int_{B_r-B_{r_0+1}} |\om|^{\frac{2m}{m-2}}\right)^\frac{m-2}{m}
\leq  \ti A\int_{B_{r_0+1}-B_{r_0}} |\om|^2 +  \frac{\ti A}{ r^2} \int_{B_{2r}-B_r} |\om|^2.
\end{equation*}
Taking $r\to \infty$ and using that $|\om|\in L^2(M)$ we have
\begin{equation}\label{first_estimates}
\left(\int_{M-B_{r_0+1}} |\om|^{\frac{2m}{m-2}}\right)^\frac{m-2}{m}
\leq \ti A\int_{B_{r_0+1}-B_{r_0}} |\om|^2.
\end{equation}
Using Holder's Inequality we obtain
\begin{eqnarray*}\label{estimate-r0}
\int_{B_{r_0+2}} |\om|^2 
\leq \vol(B_{r_0+2})^{\frac{2}{m}}
\left(\int_{B_{r_0+2}-B_{r_0+1}} |\om|^{\frac{2m}{m-2}}\right)^\frac{m-2}{m} 
+\int_{B_{r_0+1}} |\om|^2.\nonumber
\end{eqnarray*}
Define $F=\left(1+\ti A\vol(B_{r_0+2})^\frac{2}{m}\right)$. It follows from (\ref{first_estimates}) that
\begin{eqnarray}\label{desigualdade1}
\int_{B_{r_0+2}} |\om|^2  \leq F\int_{B_{r_0+1}}|\om|^2.
\end{eqnarray}

>From inequality (\ref{bks}) we have the following: 
\begin{equation}\label{desigualdade2}
|\om|\De|\om|\geq \frac{1}{m-1}|\na|\om||^2 - \Ga|\om|^2,
\end{equation}
where $\Ga:M\to [0,\infty)$ is the function given by

\begin{equation*}
\Ga=\left|(m-1)\left(|H|^2-k^2\right)-\frac{m-1}{m}|\Phi|^2-
\frac{(m-2)\sqrt{m(m-1)}}{m}|H| |\Phi|\right|.
\end{equation*}
Fix $x\in M$ and take $\zeta\in C_0^1(B_1(x))$.
Multiplying both sides of (\ref{desigualdade2}) by $\zeta^2|\om|^{p-2}$, 
with $p\geq 2$, and integrating by parts we obtain

\begin{eqnarray}\label{eta h}
-2\int_{B_1(x)} \zeta |\om|^{p-1}\lan \na\zeta,\na |\om|\ran 
&\geq& \big(p-1+\frac{1}{m-1}\big)\int_{B_1(x)} |\om|^{p-2}\zeta^2|\na |\om||^2 \nonumber \\ 
&&- \int_{B_1(x)}\Ga \zeta^2 |\om|^p.
\end{eqnarray}
Using Cauchy-Schwarz inequality (with $\varep=m-1$) we have
\begin{eqnarray*}\label{cs} 
-2\zeta |\om|^{p-1}\lan \na\zeta,\na |\om|\ran 
&=& 2\lan-|\om|^{\frac{p}{2}}\na\zeta,|\om|^{\frac{p}{2}-1}\zeta\na |\om|\ran\\
& \leq &(m-1)|\om|^p|\na\zeta|^2 + \frac{1}{m-1}|\om|^{p-2}\zeta^2|\na |\om||^2.
\end{eqnarray*} 
Applying this inequality in (\ref{eta h}) we obtain
\begin{equation}\label{pep}
(p-1) \int_{B_1(x)} |\om|^{p-2}\zeta^2|\na|\om||^2 
\leq \int_{B_1(x)}\Ga \,\zeta^2 |\om|^p +(m-1)\int_{B_1(x)} |\om|^p|\na \zeta|^2.
\end{equation}
Using Cauchy-Schwarz inequality (with $\varep=1/2$) we have
\begin{eqnarray}\label{pep2}
\int_{B_1(x)}|\na(\zeta|\om|^{\frac{p}{2}})|^2
&\leq&  \left(p+1\right)\int_{B_1(x)}|\om|^p|\na \zeta|^2 \nonumber \\
&&+\frac{p}{4}\left(p+1\right)\int_{B_1(x)}|\om|^{p-2}\zeta^2|\na |\om||^2.
\end{eqnarray}
Thus, using (\ref{pep}) and (\ref{pep2}), we obtain
\begin{eqnarray}\label{AB}
\int_{B_1(x)} |\na(\zeta|\om|^\frac{p}{2})|^2 
\leq \int_{B_1(x)} \m A\,\Ga\,\zeta^2|\om|^p + \m B\,|\om|^p|\na \zeta|^2,
\end{eqnarray}
where 
\begin{eqnarray*}
\m A&=&\frac{p}{4}\left(p+1\right)(p-1)^{-1} \leq p\quad \textrm{ and }\\
\m B &=& p+1+(m-1)\m A \leq 1+mp.
\end{eqnarray*}
In particular, $\m A,\m B \leq 2mp$, since $p\geq 2$. 
Applying the Hoffmann-Spruck Inequality \cite{hs} to the function $\zeta|\om|^{p/2}$ and using
(\ref{AB}) we get 
\begin{eqnarray*}\label{sob}
S^{-1}\left(\int_{B_1(x)} (\zeta|\om|^\frac{p}{2})^{\frac{2m}{m-2}}\,\right)^{\frac{m-2}{m}}
&\leq&\int_{B_1(x)} |\na (\zeta |\om|^{\frac{p}{2}})|^2 + 
 |H|^2 (\zeta |\om|^{\frac{p}{2}})^2\nonumber\\
&\leq& \int_{B_1(x)} \left((\m A \Ga + |H|^2)\zeta^2 + 
\m B |\na \zeta|^2\right)|\om|^p\nonumber.
\end{eqnarray*}
For simplicity, we write
\begin{eqnarray}\label{sob}
\left(\int_{B_1(x)} (\zeta|\om|^\frac{p}{2})^{\frac{2m}{m-2}}\,\right)^{\frac{m-2}{m}}
\leq 2mpS \int_{B_1(x)} \left(G\zeta^2 + |\na \zeta|^2\right)|\om|^p,  
\end{eqnarray}
where $G=\Ga +|H|^2$.

\medskip
Given $k\geq 0$ integer, we set 
$p_k=\frac{2m^k}{(m-2)^k}$ and $\rho_k=\frac{1}{2}+\frac{1}{2^{k+1}}$.
Take a function $\zeta_k\in C^{\infty}_0(B_{\rho_k}(x))$ satisfying: 
$0\leq \zeta_k\leq 1$, $\zeta_k=1$ in 
$B_{\rho_{k+1}}(x)$ and $|\na \zeta_k|\leq 2^{k+3}$. 
Using (\ref{sob}) with $p=p_k$ and $\zeta=\zeta_k$, we obtain
\begin{eqnarray*}
\left(\int_{B_{\rho_{k+1}}(x)} |\om|^{p_{k+1}}\right)^{\frac{1}{p_{k+1}}}
&\leq& (2mp_{k} S)^{\frac{1}{p_k}}
\left(\int_{B_{\rho_k}(x)} (4^{k+3}+G)|\om|^{p_k}\right)^{\frac{1}{p_k}}\nonumber\\
&\leq&\big(2mp_{k}S 4^k(4^3+\sup_{B_1(x)}G)\big)^{\frac{1}{p_k}}
\left(\int_{B_{\rho_k}(x)}|\om|^{p_k}\right)^{\frac{1}{p_k}}\nonumber \\ 
&\leq& (p_{k})^{\frac{1}{p_k}}(4^{k+k_0})^{\frac{1}{p_k}} 
\left(\int_{B_{\rho_k}(x)}|\om|^{p_k}\right)^{\frac{1}{p_k}},
\end{eqnarray*}
where $k_0$ is an integer
such that  $2m S(4^3+\sup_{B_1(x)}G)\leq 4^{k_0}$. 
By recurrence 	we obtain
\begin{eqnarray*}
||\om||_{L^{p_{k+1}}({{B_{\frac{1}{2}}}(x)})}
\leq \prod_{\ell=0}^k p_{\ell}^{\frac{1}{p_\ell}}\,4^{\frac{\ell}{p_\ell}} 4^{\frac{k_0}{p_\ell}}
||\om||_{L^2(B_1(x))}.
\end{eqnarray*}
Notice that $p_{\ell}^{\frac{1}{p_\ell}},
\,4^{\frac{\ell}{p_\ell}} \leq B^{\frac{\ell}{2}a^\ell}$ and 
$4^{\frac{k_0}{p_\ell}}\leq B^{ba^\ell}$, where $a=(m-2)/m$ and $B$, $b$ 
are suitable positive constants. Thus
$$\prod_{\ell=0}^\infty  p_{\ell}^{\frac{1}{p_\ell}}\,4^{\frac{\ell}{p_\ell}} 4^{\frac{k_0}{p_\ell}}
\leq B^{\sum_\ell  a^{\ell}(\ell+b)}<\m D,$$ where $\m D>0$ depends only on 
$m$ and $\sup_{B_1(x)}G$. Taking $k\to \infty$, we obtain 
\begin{eqnarray}\label{converge}
||\om||_{L^{\infty}({{B_{\frac{1}{2}}}(x)})}
\leq \m D\, ||\om||_{L^2(B_1(x))}.
\end{eqnarray}

Now, take  $y\in {\overline {B}}_{r_0+1}$ so that 
$\sup_{B_{r_0+1}}|\om|^2=|\om(y)|^2$.
Since $B_1(y)\subset B_{r_0+2}$, using (\ref{converge}), we obtain
\begin{equation*}
\sup_{B_{r_0+1}}|\om|^2\leq \m D\, ||\om||_{L^2(B_1(y))}^2
\leq \m D\, ||\om||_{L^2(B_{r_0+2})}^2.
\end{equation*}
Thus, from  (\ref{estimate-r0}), we have
\begin{equation}\label{final}
\sup_{B_{r_0+1}}|\om|^2 \leq \m E ||\om||_{L^2(B_{r_0+1})}^2,
\end{equation}
for all $\om\in\m H^1(M)$,
where $\m E>0$ depends on $m$, $\vol(B_{r_0+2})$ and $\sup_{B_{r_0+2}}G$.

Finally, let $\m V$ be any finite-dimensional subspace of $\m H^1(M)$. 
According to Lemma 11 of \cite{l2} there exists  $\om\in\m V$ such that

\begin{equation}\label{final-2}
\frac{\dim \m V}{\vol(B_{r_0+1})} ||\om||_{L^2(B_{r_0+1})}^2
\leq  \sup_{B_{r_0+1}}|\om|^2 \left(\min \left\{m, \dim \m V\right\}\right).
\end{equation}
Using (\ref{final}) and (\ref{final-2}), we have that $\dim \m V\leq c_0$,
where $c_0$ depends only on $m$, $\vol(B_{r_0+2})$ and 
$\sup_{B_{r_0+2}}G$. 
This implies that  $\m H^1(M)$  has finite dimension. 
This concludes the proof of Theorem \ref{t1}.

\medskip


\section*{Acknowledgement} The first author wishes to thank Haifeng Xu for his interesting comments
about this manuscript. The second and the third authors thank Walcy Santos and 
Maria Fernanda Elbert for helpful suggestions during the preparation of this article.

\vspace{0.7cm}

\begin{small}

\begin{tabular}{l}
Marcos Petr\'ucio Cavalcante\\
Universidade Federal de Alagoas\\
Instituto de Matem\'{a}tica\\
57072-900 Macei\'o-AL\\
Brazil\\
\verb+marcos.petrucio@pq.cnpq.br+\\
\end{tabular}\

\vspace{1cm}

\begin{tabular}{l}
Heudson Mirandola\\
Universidade Federal do Rio de Janeiro\\
Instituto de Matem\'{a}tica\\
21945-970 Rio de Janeiro-RJ\\
Brazil\\
\verb+mirandola@im.ufrj.br+\\
\end{tabular}\

\vspace{1cm}

\begin{tabular}{l}
Feliciano Vit\'orio\\
Universidade Federal de Alagoas\\
Instituto de Matem\'{a}tica\\
57072-900 Macei\'o-AL\\
Brazil\\
\verb+feliciano.vitorio@pq.cnpq.br+\\
\end{tabular}\

\end{small}


\begin{thebibliography}{99}
\bibitem{a} Anderson, M., 
{\it The compactification of a minimal submanifold in Euclidean space by the Gauss map}, 
Preprint I.H.ƒ.S, 1985.

\bibitem{a2}
Anderson, M., 
{\it $L^2$ harmonic forms and a conjecture of Dodziuk-Singer}. 
Bull. Amer. Math. Soc. (N.S.) 13 (1985), no. 2, 163 -- 165.



%


\bibitem{bc}  Bessa, G. P. and Costa, M. S., 
{\it On submanifolds with tamed second fundamental form}. 
Glasg. Math. J. 51 (2009), no. 3, 669 -- 680.

\bibitem{BdCS}  B\'erard, P.,  do Carmo, M. and Santos, W.,
{\it Complete hypersurfaces with constant mean curvature and finite total curvature.}
Ann. Global Anal. Geom. 16 (1998), no. 3, 273Ð-290. 

\bibitem{b} Bochner. S.,
{\it Vector fields and Ricci curvature}. 
Bull. Amer. Math. Soc.  52 (1946), 776 -- 797.

\bibitem{bmj} Bessa, G. P., Jorge, L. and Montenegro, J. F.,
{\it Complete submanifolds of $\real^n$ with finite topology}.  
Comm. Anal. Geom. 15 (2007), no. 4, 725 -- 732. 



\bibitem{cwx} do Carmo, M. P., Wang, Q. and Xia, C., 
{\it Complete submanifolds with bounded mean curvature in a Hadamard manifold}.
 J. Geom. Phys. 60 (2010), no. 1, 142 -- 154.


\bibitem{dccs} do Carmo, M. P., Cheung, L-F. and Santos, Walcy
{\it On the compactness of constant mean curvature hypersurfaces with finite total curvature}.
Arch. Math. (Basel) 73 (1999), no. 3, 216--222. 

\bibitem{c} Carron, G., 
{\it $L^2$-Cohomologie et in\'egalit\'es de Sobolev}. 
Math. Ann. 314 (1999), 613 -- 639.


%
%





\bibitem{co} Chern, S. S. and Osserman, R., 
{\it Complete minimal surfaces in euclidean n space}. 
J. Analyse Math. 19 (1967) 15 -- 34.









\bibitem{fl} Fu, H-P and  Li, Z-Q., 
{\it $L^2$ harmonic $1$-forms on complete submanifolds in Euclidean space}. 
Kodai Math. J. 32 (2009), 432 -- 441.


\bibitem{fx} Fu, H-P. and Xu, H-W., 
{\it Total curvature and $L^2$ harmonic $1$-forms on complete submanifolds in space forms}. 
Geom. Dedicata  144  (2010), 129 -- 140.




\bibitem{h} Huber, A., 
{\it On subharmonic functions and differential geometry in the large}. 
Comment. Math. Helv., 52 (1957), 13 -- 72 .


\bibitem{hs} Hoffman, D. and Spruck, J., 
{\it Sobolev and isoperimetric inequalities for Riemannian submanifolds}. 
Comm. Pure. Appl. Math. 27 (1974), 715 -- 727.



 
\bibitem{l2} Li, P., 
{\it On the Sobolev constant and the $p$-spectrum of a compact Riemannian manifold}. 
Ann. Sci. \'Ecole Norm. Sup. (4) 13 (1980), no. 4, 451 -- 468. 

\bibitem{li} Li, P., 
{\it Lecture notes on geometric analysis}. 
Lecture Notes Series, 6. Seoul National University, Research Institute of Mathematics, 
Global Analysis Research Center, Seoul, 1993.




\bibitem{lt} Li, P. and Tam, L. F., 
{\it Harmonic functions and the structure of complete manifolds}. 
J. Diff. Geom. 35 (1992), no. 2, 359 -- 383.


\bibitem{lw} Li, P. and Wang, J.P., 
{\it Minimal hypersurfaces with finite index}.  
Math. Res. Lett. 9 (2002), no. 1, 95 -- 103.


\bibitem{ms} M\"uller, S. and \u Sver\'ak, V., 
{\it On surfaces of finite total curvature}. 
J. Differential Geom. 42 (1995), no. 2, 229 -- 258.



\bibitem{ni}  Ni, L., 
{\it Gap theorems for minimal submanifolds in $\real^{n+1}$}. 
Comm. Anal. Geom. 9 (2001), 641 -- 656.


\bibitem{do}  Oliveira Filho, G., 
{\it Compactification of minimal submanifolds of hyperbolic space}.
Comm. Anal. Geom. 1 (1993), no. 1, 1 -- 29. 


\bibitem{o} Osserman, R., 
{\it Global properties of minimal surfaces in $E^3$ and $E^n$}.
 Ann. of Math., 80  (1964), 340--364.


\bibitem{seo} Seo, K., 
{\it Minimal submanifolds with small total scalar curvature in Euclidean space}. 
Kodai Math. J., 31 (2008), 113 -- 119.


\bibitem{seo2} Seo, K., 
{\it Rigidity of minimal submanifolds in hyperbolic space}. 
Arch. Math. 94 (2010), 173 -- 181.


\bibitem{sx} Shiohama, K. and Xu, H., 
{\it The topological sphere theorem for complete submanifolds}. 
Compositio Math. 107 (1997), 221 -- 232.


\bibitem{w} Wang, X., {\it On conformally compact Einstein manifolds}. 
Math. Res. Lett.  8  (2001),  no. 5-6, 671 -- 688. 


\bibitem{wh} White, B., {\it Complete surfaces of finite total curvature}. J. Diff. Geom. 26 (1987), 315--326.

\end{thebibliography}
\end{document}